\documentclass[a4paper]{amsart}

\usepackage{amssymb}
\usepackage{amsthm}
\usepackage{amscd}
\usepackage{amsmath}
\usepackage[all]{xy}
\newtheorem{theorem}{Theorem}[subsection]
\newtheorem{claim}[theorem]{Claim}
\newtheorem{lemma}[theorem]{Lemma}
\newtheorem{proposition}[theorem]{Proposition}
\newtheorem{cor}[theorem]{Corollary}

\theoremstyle{definition}
\newtheorem{definition}[theorem]{Definition}
\newtheorem{remark}[theorem]{Remark}

\newcommand{\Z}{{\mathbb Z}}
\newcommand{\ZZ}{{\mathbb Z}}
\newcommand{\C}{{\mathbb C}}
\newcommand{\CC}{{\mathbb C}}

\newcommand{\IR}{{\mathbb R}}

\newcommand{\I}{\mathbf{1}}

\newcommand{\lk}{\ell\mathit{k}}
\newcommand{\bord}{\partial}

\newcommand{\SL}{{\mathrm{SL}_2(\C)}}
\newcommand{\PSL}{{\mathrm{PSL}_2(\C)}}

\newcommand{\T}{\mathcal{T}}

\newcommand{\sll}{\mathfrak{sl}_2(\CC)}
\newcommand{\tsll}{\widetilde{\mathfrak{sl}}_2(\CC)}

\newcommand{\trace}{{\rm Tr}\,}
\newfont{\bg}{cmr10 scaled\magstep4}

\newcommand{\im}{\mathop{\mathrm{im}}\nolimits}

\newenvironment{proofclaim}{%
  \begin{proof}%
}{%
  \end{proof}%
}

\begin{document}

\thanks{This research is partially supported by 
the 21st century COE program at Graduate School of Mathematical
Sciences, the University of Tokyo}

\title{
Limit values of the non-acyclic Reidemeister torsion 
for knots}


\author{Yoshikazu Yamaguchi}
\address{Graduate School of Mathematical Sciences, 
University of Tokyo, 3-8-1 Komaba Meguro, 
Tokyo 153-8914, Japan}
\email{shouji@ms.u-tokyo.ac.jp}

\date{}

\begin{abstract}
We consider the Reidemeister torsion associated with $\SL$-representations
of a knot group.
A bifurcation point in the $\SL$-character variety of a knot group
is a character which is given by both 
an abelian $\SL$-representation and 
a non-abelian one.
We show that there exist limits of the non-acyclic Reidemeister torsion 
at bifurcation points and the limits are expressed by using the derivation of 
the Alexander polynomial of the knot in this paper.
\end{abstract}

\maketitle



\section{Introduction}
The Reidemeister torsion is an invariant of 
a CW-complex and a representation of its fundamental group.
For a knot exterior and an abelian representation, 
the Reidemeister torsion is essentially equal to the Alexander polynomial,
see Milnor~\cite{Milnor:duality, Milnor:covering} and Turaev~\cite{Turaev:2001}.
In the case of a non-abelian representation, 
the Reidemeister torsion is related to the theory of the twisted Alexander invariant,
see Kirk and Livingston~\cite{ KL}, Kitano~\cite{Kitano}, Lin~\cite{Lin} and Wada~\cite{wada}.

The Reidemeister torsion is invariant 
under taking conjugation of a representation. 
And in the case of knot exteriors the Reidemeister torsion may be regarded as a function on a space corresponding to a suitable quotient of the $\SL$-representations of the knot group by conjugation, we can find that this point of view was introduced in Porti~\cite{Porti:1997}.
Following Morgan and Shalen~\cite{MS},
we consider the $\SL$-character variety of the knot group as a suitable quotient. 
In general, the $\SL$-character variety of a knot group has 
many components.
These components are roughly classified into two types.
One consists of the characters of abelian representations.
The other consists of the characters of non-abelian representations. 
We respectively call these sets the {\it abelian part} and 
the {\it non-abelian part} 
of the character variety.
It is known that the abelian part intersects with the non-abelian part.
These intersection points are called {\it bifurcation points}.
The purpose of this paper is to show that 
the Reidemeister torsion of 
non-abelian representations is given  by using the Alexander 
polynomials at a bifurcation point as follows.

Let $K$ be a knot in a homology three sphere.
A bifurcation point of the $\SL$-character variety of $K$
corresponds to a root of the Alexander polynomial, see Burde~\cite{BurdeBifurcation} and Klassen~\cite{Klassen}.
In particular, the bifurcation point corresponding to a simple root of the Alexander polynomial is a smooth point of 
the $\SL$-character variety (see Heusener, Porti and Su\'arez~\cite{HPS}).
We can construct a function on each of the abelian and 
non-abelian part of the character variety by using the Reidemeister torsion.
The function on the abelian part is given 
by the Reidemeister torsion for abelian representations. 
In fact, this function is expressed by using the Alexander polynomial of $K$
and it has zeros at bifurcation points
(see Milnor~\cite{Milnor:duality, Milnor:covering} and Turaev~\cite{Turaev:2001}).
The other function on the non-abelian part is given by 
the {\it non-acyclic} Reidemeister torsion for non-abelian representations,  Dubois~\cite{JDFourier, JDFibre}, Porti~\cite{Porti:1997} and Yamaguchi~\cite{Yamaguchi} deal with Reidemeister torsion in such a light.
Though the function on the non-abelian part is partially defined and it is not defined on bifurcation points,
we can consider limits of the non-acyclic Reidemeister torsion at 
bifurcation points.

We will show that 
if a bifurcation point corresponds to a simple root of 
the Alexander polynomial of $K$, then 
there exists 
the limit of the non-acyclic Reidemeister torsion at the bifurcation point
and 
its limit is expressed as
the differential coefficient of the function defined on the abelian part 
at this point
(Theorem \ref{thm:Main_theorem}).

This fact had been first conjectured by Dubois and Kashaev.
The author proved it for a knot in $S^3$ at first.
Dubois pointed out that the proof may be extended to a knot 
in a homology three sphere.
This theorem is applied in the paper of 
Dubois and Kashaev~\cite{DuboisKashaev}.


This paper is organized as follows.
In Section \ref{section:preliminary}, we recall the needed notions
of the $\SL$-character variety of a knot group and 
the Reidemeister torsion for knot exteriors.
In Section \ref{section:torsion_at_bifurcation},
we prove that limits of
the non-acyclic Reidemeister torsion of a knot exterior at bifurcation points
are obtained from the derivation of the Alexander polynomial of the knot.
We discuss the existences of limits of the non-acyclic Reidemeister torsion
in Subsection \ref{section:existence_limit}.
We give a formula of these limits 
in Subsection \ref{section:limit_torsion_formula}.
This formula implies that some property,
called {\it $\lambda$-regularity},
 which holds 
on irreducible characters near a bifurcation point
can be extended to the bifurcation point. 
This is shown in Section \ref{section:reference_generators}.

\subsection*{Acknowledgements}
The author would like to express sincere gratitude to Mikio Furuta 
for his suggestions and helpful discussions.
He is thankful to Hiroshi Goda, Takayuki Morifuji, Teruaki Kitano, 
Masaaki Suzuki and Yuya Koda
for helpful suggestions.
The author would like to thank 
J\'er\^ome Dubois for his helpful advices.
He also would like to thank the referee for his/her careful reading 
and appropriate advices.

\section{Preliminary}
\label{section:preliminary}

\subsection{Review on bifurcation points}
\label{section:review_bifurcation}
Let $K$ be a knot in a homology three sphere $M$, $M_K$ its exterior
and $R(\pi_1(M_K), \SL)$ denote
the set of $\SL$-representations of $\pi_1(M_K)$.

A representation $\rho$ is called {\it abelian\/}
if its image $\rho(\pi_1(M_K))$ is an abelian subgroup of $\SL$.
A representation $\rho$ is called {\it reducible\/}
if there exists a proper subspace $U$ of $\C^2$ such that
$\rho(\gamma)(U) \subset U$ for any $\gamma \in \pi_1(M_K)$.
A representation $\rho$ is called {\it irreducible\/} if it is not reducible.
We let $R^{irr}(\pi_1(M_K), \SL)$ denote the set of irreducible ones.
Note that all abelian representations are reducible but the converse is false in general.

Associated to the representation $\rho \in R(\pi_1(M_K), \SL)$
is its {\it character} a map $\chi_{\rho}$ from $\pi_1(M_K)$ into $\C$, 
defined by $\chi_{\rho}(\gamma)=\trace(\rho(\gamma))$.
Following Morgan and Shalen~\cite{MS},
we will focus on the {\it character variety\/} 
which is the set of ${\it characters\/}$ of 
$\SL$-representations of $\pi_1(M_K)$.
Let $X(M_K)$ denote the character variety of $\pi_1(M_K)$.
In some sense, $X(M_K)$ is 
the \lq\lq algebraic quotient\rq\rq\, of $R(\pi_1(M_K), \SL)$ 
by $\PSL$ because the quotient $R(\pi_1(M_K), \SL)/\PSL$ is not Hausdorff in general.
We let $\pi$ denote the projection, $R(\pi_1(M_K), \SL) \to X(M_K)$, defined 
by $\rho \mapsto \chi_{\rho}$.
It is known that 
$R(\pi_1(M_K), \SL)$ and $X(M_K)$ 
have the structure of complex algebraic affine sets 
and for each $\gamma \in \pi_1(M_K)$ the function 
$I_{\gamma}:X(M_K) \to \C, \chi_\rho \mapsto \trace(\rho(\gamma))$
is a regular function.
Two irreducible representations of $\pi_1(M_K)$ with 
the same character are conjugate by an element of $\SL$ 
(see Culler and Shalen~\cite[Proposition 1.5.2]{CullerShalen}).
Let $X^{irr}(M_K)$ denote $\pi(R^{irr}(\pi_1(M_K), \SL))$.
The subsets $R^{irr}(\pi_1(M_K), \SL) \subset R(\pi_1(M_K), \SL)$
and $X^{irr}(M_K) \subset X(M_K)$ are Zariski-open.
(For the details, see Morgan and Shalen~\cite{MS}.)


The character variety $X(M_K)$ has several components.
Let $X^{ab}(M_K)$ be 
the image under $\pi$ of the subset of abelian $\SL$-representations of $\pi_1(M_K)$ and $X^{nab}$ the image of the subset of non-abelian ones.
We call $X^{ab}(M_K)$ (resp. $X^{nab}(M_K)$) 
the {\it abelian} (resp. {\it non-abelian}) part of $X(M_K)$.
\begin{definition}
If there exist intersection points between 
the abelian part $X^{ab}(M_K)$ and the non-abelian part $X^{nab}(M_K)$
in $X(M_K)$,
then these intersection points are called {\it bifurcation points}. 
\end{definition}

It is well known that 
$\pi_1(M_K) / [\pi_1(M_K),\, \pi_1(M_K)] 
\cong H_1(M_K;{\mathbb Z}) 
\cong {\mathbb Z}$
is generated by the meridian $\mu$ of $K$.  
\begin{remark}\label{rem:maximal_sub}
In $\SL$ there exist, up to conjugation, 
only two maximal abelian subgroups Hyp$(=$hyperbolic$)$
and Para$(=$parabolic$)$; 
they are given by
\begin{align*}
{\rm Hyp} 
&:=
\left\{\left.
\left(
  \begin{array}{cc}
     c  & 0 \\
     0  & c^{-1}
  \end{array}
\right) \in \SL
\,\right|\,
c \in \C^* = \C \setminus \{0\}
\right\},\\
{\rm Para} 
&:=
\left\{\left.
\pm
\left(
  \begin{array}{cc}
     1 & \omega \\
     0 & 1
  \end{array}
\right) \in \SL
\,\right|\,
\omega \in \C
\right\}.
\end{align*}
\end{remark}

As a consequence, 
each abelian representation of $\pi_1(M_K)$ in $\SL$
is conjugate either to
\[
\varphi_z :
  \pi_1(M_K) \ni \mu \mapsto 
  \left(
    \begin{array}{cc}
      e^{z} & 0 \\
      0     & e^{-z}    
    \end{array}
  \right)
  \in \SL
\]
with $z \in \C$ if it is hyperbolic,
or
to a representation $\rho$ with
$\displaystyle{
\rho(\mu)
=
\pm
\left(
  \begin{array}{cc}
    1 & 1 \\
    0 & 1
  \end{array}
\right)
}$ if it is parabolic.

The non-abelian part $X^{nab}(M_K)$ includes 
the irreducible characters $X^{irr}(M_K)$.
It is known that an element of $X^{irr}(M_K)$ is a smooth point in 
the complex affine variety $X(M_K)$, for example see Porti~\cite[Proposition 3.5]{Porti:1997}.
We focus on the bifurcation points which are limits of paths in $X^{irr}(M_K)$.
Such bifurcation points are related to 
roots of the Alexander polynomial $\Delta_K(t)$ of $K$.
This is a well-known result of
Burde~\cite{BurdeBifurcation}, de Rham~\cite{deRham}
if $K$ is a knot in $S^3$.

\begin{lemma}[Corollary 4.3 in Heusener--Porti--Su\'arez~\cite{HPS}, Klassen~\cite{Klassen}]
\label{lemma:bifurcation_point}
Let $z_0$ be a complex number.
There is a reducible non-abelian representation $\rho_{z_0}$
such that $\chi_{\rho_{z_0}} = \chi_{\varphi_{{z_0}}}$
if and only if  $\Delta_K(e^{2z_0})=0$.
\end{lemma}

It is also known that the following theorem holds.

\begin{theorem}[Theorem 1.1 in Heusener--Porti--Su\'arez~\cite{HPS}]
\label{thm:bifurcation_point}
Let $z_0$ be a complex number such that $\Delta_K(e^{2z_0})=0$ and 
$\rho_{z_0}$ a reducible non-abelian representation 
such that $\chi_{\rho_{z_0}} = \chi_{\varphi_{{z_0}}}$.
If $e^{2z_0}$ is a simple root of $\Delta_K(t)$, 
then the representation 
$\rho_{z_0}$ is the limit of a sequence of irreducible ones.
More precisely, $\rho_{z_0}$ is a smooth point of the $\SL$-representation 
variety of $\pi_1(M_K)$;
it is contained in a unique irreducible four-dimensional component
of the $\SL$-representation variety. 
\end{theorem}

Heusener, Porti and Su\'arez also showed that the character of $\rho_{z_0}$ is a smooth point of the 
$\SL$-character variety $X(M_K)$ (see Theorem 1.2 in Heusener--Porti--Su\'arez~\cite{HPS}).

We will consider bifurcation points corresponding to 
simple roots of the Alexander polynomial $\Delta_K(t)$.
These bifurcation points are limits of paths in $X^{irr}(M_K)$.

\subsection{
Review on the Reidemeister torsion}
\label{Review_Reidemeister_torsion}
\subsubsection*{Torsion of a chain complex}
Let 
$C_* = 
(0 \to 
 C_n \xrightarrow{d_n} 
 C_{n-1} \xrightarrow{d_{n-1}} 
 \cdots \xrightarrow{d_1} 
 C_0 \to 0)$ 
be a chain complex of finite dimensional vector spaces over $\C$. 
Choose a basis $\mathbf{c}^i$ for $C_i$ and  a basis $\mathbf{h}^i$ for 
the $i$-th homology group $H_i = H_i(C_*)$. 
The torsion of $C_*$ with respect to these choice of bases is defined as follows.

Let $\mathbf{b}^i$ be a sequence of vectors in $C_{i}$ 
such that $d_{i}(\mathbf{b}^i)$ is a basis of 
$B_{i-1}= \im(d_{i} \colon C_{i} \to C_{i-1})$ and 
let $\widetilde{\mathbf{h}}^i$ denote a lift of $\mathbf{h}^i$ in 
$Z_i = \ker(d_{i} \colon C_i \to C_{i-1})$. 
The set of vectors 
$d_{i+1}(\mathbf{b}^{i+1})\widetilde{\mathbf{h}}^i\mathbf{b}^i$ 
is a basis of $C_i$. 
Let 
$[d_{i+1}(\mathbf{b}^{i+1})\widetilde{\mathbf{h}}^i\mathbf{b}^i/\mathbf{c}^i] 
\in \C^*$ 
denote the determinant of the transition matrix 
between those bases 
(the entries of this matrix are coordinates of vectors in 
$d_{i+1}(\mathbf{b}^{i+1})\widetilde{\mathbf{h}}^i\mathbf{b}^i$ 
with respect to $\mathbf{c}^i$). 
The \emph{sign-determined Reidemeister torsion} of $C_*$ 
(with respect to the bases $\mathbf{c}^*$ and $\mathbf{h}^*$) 
is the following alternating product (see Turaev~\cite[Definition 3.1]{Turaev:2001}):
\begin{equation}
\label{Def:RTorsion}
\mathrm{Tor}(C_*, \mathbf{c}^*, \mathbf{h}^*) 
= (-1)^{|C_*|} \cdot 
\prod_{i=0}^n 
[
d_{i+1}(\mathbf{b}^{i+1})\widetilde{\mathbf{h}}^i\mathbf{b}^i/\mathbf{c}^i
]^{(-1)^{i+1}} 
\in \C^*.
\end{equation}
Here  $$|C_*| = \sum_{k\geqslant 0} \alpha_k(C_*) \beta_k(C_*),$$ 
where 
$\alpha_i(C_*) = \sum_{k=0}^i \dim C_k$,  $\beta_i(C_*)  = \sum_{k=0}^i \dim H_k$.

The torsion $\mathrm{Tor}(C_*, \mathbf{c}^*, \mathbf{h}^*)$ 
does not depend on the choices of $\mathbf{b}^i$ and $\widetilde{\mathbf{h}}^i$.
Further observe that if $C_*$ is acyclic (i.e., if $H_i = 0$ for all $i$), 
then $|C_*| = 0$.

\subsubsection*{Torsion of a CW-complex}
Let $W$ be a finite CW-complex, $V$ a finite dimensional vector space
over $\C$
and $\rho$ a homomorphism from $\pi_1(W)$ to $Aut(V)$.
We define the local system of $W$ to be
\[
C_*(W; V_\rho) =  V_\rho \otimes_{\ZZ[\pi_1(W)]} C_*(\widetilde{W}; \ZZ).
\]
Here $C_*(\widetilde{W}; \ZZ)$ is the complex of 
the universal cover $\widetilde{W}$ with integer coefficients. 
This space is in fact a left $\ZZ[\pi_1(W)]$-module 
(via the action of $\pi_1(W)$ on $\widetilde{W}$ as the covering group). 
And $V_\rho$ denotes the right $\ZZ[\pi_1(W)]$-module 
via the homomorphism $\rho$, 
i.e., the action is given by
$v \cdot \gamma = \rho(\gamma)^{-1}(v)$ 
for any $v \in V$ and $\gamma \in \pi_1(W)$.
This chain complex $C_*(W; V_\rho)$ computes the homology of the local system. 
We let $H_*(W; V_\rho)$ denote this homology.

Let $\{e^{(i)}_1, \ldots, e^{(i)}_{n_i}\}$ be 
the set of $i$-dimensional cells of $W$. 
We lift them to the universal cover and 
we choose an arbitrary order and an arbitrary orientation for 
the cells $\left\{ {\tilde{e}^{(i)}_1, \ldots, \tilde{e}^{(i)}_{n_i}} \right\}$. 
If $\mathcal{B} = \{\mathbf{f_1}, \ldots, \mathbf{f_m}\}$ is 
an orthonormal basis of $V$, where $m$ is the dimension of $V$, 
then we consider the corresponding basis over $\CC$
$$\mathbf{c}^{i}_{\mathcal{B}} = 
\left\{
  \mathbf{f_1} \otimes \tilde{e}^{(i)}_{1}, 
  \ldots, 
  \mathbf{f_m} \otimes \tilde{e}^{(i)}_{1}, \ldots, 
  \mathbf{f_1} \otimes \tilde{e}^{(i)}_{n_i}, 
  \ldots, 
  \mathbf{f_m} \otimes \tilde{e}^{(i)}_{n_i}
\right\}$$ 
of 
$C_i(W; V_\rho)$. 
Now choosing for each $i$ 
a basis $\mathbf{h}^{i}$ for the homology group $H_i(W; V_\rho)$, 
we can compute 
\[
\mathrm{Tor}(C_*(W; V_\rho), \mathbf{c}^*_{\mathcal{B}}, \mathbf{h}^{*}) \in \C^*.
\]

The cells 
$\{ \tilde{e}^{(i)}_j \}^{}_{0 \leqslant i \leqslant \dim W, 1 \leqslant j \leqslant n_i}$
 are in one-to-one correspondence with the cells of $W$, 
their order and orientation induce an order and an orientation for the cells 
$\{ e^{(i)}_j \}^{}_{i, j}$, 
where $0 \leqslant i \leqslant \dim W$ and $1 \leqslant j \leqslant n_i$. 
Again, corresponding to these choices, 
we get a basis $c^i$ over $\IR$ for $C_i(W; \IR)$. 

Choose an \emph{homology orientation} of $W$, 
which is an orientation of the real vector space 
$H_*(W; \IR) = \bigoplus_{i\geqslant 0} H_i(W; \IR)$. 
Let $\mathfrak{o}$ denote this chosen orientation. 
Provide each vector space $H_i(W; \IR)$ with a reference basis $h^i$ 
such that the basis 
$\left\{ {h^0, \ldots, h^{\dim W}} \right\}$ of $H_*(W; \IR)$ is 
{positively oriented} with respect to $\mathfrak{o}$. 
We set 
\[
\tau_0 = 
\mathrm{sgn}\left(\mathrm{Tor}(C_*(W; \IR), c^*, h^{*})\right) \in \{\pm 1\}.
\] 

We define 
the sign--determined Reidemeister torsion for $(W, V_\rho)$ 
with respect to the homology basis $\mathbf{h}^*$ and to the homology orientation $\mathfrak{o}$
to be
\begin{equation}\label{EQ:TorsionRaff}
\mathrm{TOR}(W; V_\rho, \mathbf{h}^{*}, \mathfrak{o}) 
= \tau_0 \cdot 
\mathrm{Tor}(C_*(W; V_\rho), \mathbf{c}^*_{\mathcal{B}}, \mathbf{h}^{*}) 
\in \C^*.
\end{equation}
This definition only depends on the combinatorial class of $W$, 
the conjugacy class of $\rho$, the choice of $\mathbf{h}^{*}$ and 
the homology orientation $\mathfrak{o}$. 
It is independent of the orthonormal basis $\mathcal{B}$ of $V$, 
of the choice of the lifts $\tilde{e}^{(i)}_j$, and 
of the choice of the positively oriented basis of $H_*(W; \IR)$. 
Moreover, it is independent of the order and the orientation of the cells 
(because they appear twice). 

\begin{remark}
If the Euler characteristic of $W$ is zero, 
then
we can use any basis of $V$ 
in order to define $\mathrm{TOR}(W;V_\rho, \mathbf{h}^{*}, \mathfrak{o})$.
\end{remark}
One can prove that $\mathrm{TOR}$ is invariant under cellular subdivision, 
homeomorphism and simple homotopy equivalences. 
In fact, 
all these important invariance properties hold with the sign $(-1)^{|C_*|}$ 
in~$(\ref{Def:RTorsion})$, 
for details see Farber and Turaev~\cite[Lemma 3.3]{FarberTuraev:Poincare}.

\subsection{Review on the non-acyclic 
Reidemeister torsion for knot exteriors}
\label{section:non_acyclic_torsion}
This subsection is devoted to a detailed review of 
the constructions of the non-acyclic Reidemeister torsion 
which were made in Dubois~\cite{JDFourier} and Porti~\cite{Porti:1997}. 

Let $K$ be a knot in a homology three sphere $M$ and $M_K$ its exterior.
We let $\rho$ denote an $\SL$-representation of $\pi_1(M_K)$ and 
$Ad$ be the adjoint action of $\SL$, i.e.,
$Ad:\SL \to Aut(\sll),\, A \mapsto (Ad_A: x \mapsto AxA^{-1})$.

We define the local system $C_*(M_K; \sll_\rho)$ by 
\[
 C_*(M_K; \sll_\rho)
:= \sll_\rho \otimes_{\Z[\pi_1(M_K)]} C_*(\widetilde M_K ;\Z) 
\]
where $\widetilde M_K$ is the universal cover of $M_K$ 
and $\sll_\rho$ is the right $\Z[\pi_1(M_K)]$-module 
via the composition $Ad \circ \rho$, i.e.,
$v \cdot \gamma = Ad_{\rho(\gamma)^{-1}}(v)$
for any $v \in \sll$ and $\gamma \in \pi_1(M_K)$.
We call this local system {\it the $\sll_\rho$-twisted chain complex} of $M_K$.

We let $H_*(M_K; \sll_\rho)$ denote the homology of this local system.
It is known that $\dim_\C H_1(M_K; \sll_\rho)$ is equal to the dimension of the 
component of $X(M_K)$ which contains $\chi_{\rho}$
if $\rho$ is irreducible.
In particular, for an irreducible representation $\rho$,
$C_*(M_K;\sll_{\rho})$ is not acyclic 
since there are no $0$-dimensional components of $X(M_K)$
(see Cooper--Culler--Gillet--Long--Shalen~\cite[Proposition 2.4]{CCGLS}). 

\subsubsection*{Canonical homology orientation of knot exteriors}
We provide 
the exterior of $K$ with its \emph{canonical homology orientation} 
defined as follows (see Turaev~\cite[Section V.3]{Turaev:2002}).
We have 
$$H_*(M_K; \IR) = H_0(M_K; \IR) \oplus H_1(M_K; \IR)$$ 
and we base this $\IR$-vector space with 
$\{ \lbrack pt \rbrack, \lbrack \mu \rbrack\}$. 
Here $\lbrack pt \rbrack$ is 
the homology class of a point, 
and $\lbrack \mu \rbrack$ is 
the homology class of the meridian $\mu$ of $K$. 
This reference basis of $H_*(M_K; \IR)$ induces 
the so--called canonical homology orientation of $M_K$. 
We let $\mathfrak{o}$ denote the canonical homology orientation of $M_K$.

\subsubsection*{Regularity for representations}

In this subsection we briefly review two notions of regularity 
(see Dubois~\cite{JDFibre} and Porti~\cite{Porti:1997}). 
Let $K \subset M$ denote an oriented knot. 

 The meridian $\mu$ of $K$ is supposed to be oriented according to 
the rule $\lk(K, \mu) = +1$, 
while the preferred longitude $\lambda$ is oriented according to 
the condition $\mathrm{int}(\mu, \lambda) = +1$. 
Here $\mathrm{int}(\cdot, \cdot)$ denotes the intersection form on $\bord M_K$. 

We say that 
$\rho \in R^\mathrm{irr}(\pi_1(M_K), \SL)$ is \emph{regular} 
if $\dim_{\C} H_1(M_K; \sll_\rho) = 1$. 
This notion is invariant by conjugation and 
thus it is well-defined for irreducible characters.
Note that for a regular representation $\rho$, 
we have 
$$
\dim_{\C} H_1(M_K; \sll_\rho) = 1, \; \dim_{\C} H_2(M_K; \sll_\rho) = 1 \text{ and } 
H_j(M_K; \sll_\rho) = 0 
$$ 
for all $j \ne 1, 2$ by Porti~\cite[Corollary 3.23]{Porti:1997}.
Let $\gamma$ be a simple closed unoriented curve in $\bord M_K$. 
Among irreducible representations we focus on the $\gamma$-regular ones. 
We say that a regular representation 
$\rho : \pi_1(M_K) \to \SL$ is \emph{$\gamma$-regular}
(see Porti~\cite[Definition 3.21]{Porti:1997}), 
if
\begin{enumerate}
  \item the inclusion $\iota \colon \gamma \hookrightarrow M_K$ induces 
        a \emph{surjective} map 
        $$\iota_* \colon H_1(\gamma; \sll_\rho) 
           \to H_1(M_K; \sll_\rho);$$
  \item if $\trace(\rho(\pi_1(\bord M_K))) \subset \{\pm 2\}$, 
        then $\rho(\gamma) \ne \pm \I$.
\end{enumerate} 
It is easy to see that this notion is invariant by conjugation. 
For $\chi \in X^{\mathrm{irr}}(M_K)$ 
the notion of $\gamma$-regularity is well-defined. 

\subsubsection*{How to construct natural bases for the twisted homology}

Let $\rho$ be a regular $\SL$-representation of $\pi_1(M_K)$ 
and fix a generator $P^\rho$ of $H_0(\bord M_K; \sll_\rho)$ 
(i.e., the vector $P^\rho$ in $\sll$ satisfies the condition that 
$Ad_{\rho(g)}(P^\rho) = P^\rho$ for all $g \in \pi_1(\bord M_K)$).  

Suppose that $M$ is oriented. The exterior of a knot is thus oriented and we know that it is bounded by a $2$-dimensional torus. 
This boundary inherits an orientation by the convention 
{\it \lq\lq the inward pointing normal vector in the last position\rq\rq}. 
The usual inclusion $i\colon \bord M_K \to M_K$ induces 
(see Dubois~\cite[Lemma 5.2]{JDFourier}) 
an isomorphism $i_*\colon H_2(\bord M_K; \sll_\rho) \to H_2(M_K; \sll_\rho)$. 
Moreover, 
one can prove that 
$H_2(\bord M_K; \sll_\rho) \cong H_2(\bord M_K; \ZZ) \otimes \CC$ 
(see Dubois~\cite[Lemma 5.1]{JDFourier}). 
More precisely, 
let 
$\lbrack \bord M_K \rbrack \in H_2(\bord M_K; \ZZ)$ 
be the fundamental class induced by the orientation of $\bord M_K$, 
we have that 
$H_2(\bord M_K; \sll_\rho) 
= \CC[P^\rho \otimes \widetilde{\bord M_K}]$.

 The \emph{reference generator} of $H_2(M_K; \sll_\rho)$ is defined by 
\begin{equation}\label{EQ:Defh2}
h_{(2)}^\rho 
= i_*([P^\rho \otimes \widetilde{\bord M_K}]).
\end{equation}

Let $\rho$ be a $\lambda$-regular representation of $\pi_1(M_K)$. 
Then 
the \emph{reference generator} of $H_1(M_K; \sll_\rho)$ is defined by
\begin{equation}\label{EQ:Defh1}
h_{(1)}^\rho(\lambda) = \iota_*([P^\rho \otimes \widetilde{\lambda}]).
\end{equation}

\begin{remark}
The generator $h_{(1)}^\rho(\lambda)$ of $H_1(M_K; \sll_\rho)$ depends on 
the orientation of $\lambda$. 
If we change the orientation of the longitude $\lambda$ 
in Equation~(\ref{EQ:Defh1}),  
then the generator changes into its reverse.
\end{remark}

\begin{remark}
Note that 
$H_i(M_K; \sll_\rho)$ is isomorphic to 
the dual space of the twisted cohomology $H^i(M_K; \sll_\rho)$. 
The reference elements defined 
in Equations~(\ref{EQ:Defh2}) and (\ref{EQ:Defh1}) are dual 
from ones defined in Dubois~\cite[\S~3.4]{JDFibre}.
\end{remark}

\subsubsection*{The non-acyclic Reidemeister torsion for  knot exteriors}

Let $\rho \colon \pi_1(M_K) \to \SL$ be a $\lambda$-regular representation. 
The \emph{Reidemeister torsion $\mathbb{T}^K_\lambda$} at $\rho$ is  
defined to be 
\begin{equation}\label{Tordef}
\mathbb{T}^K_\lambda(\rho) 
= \mathrm{TOR}
\left( 
{M_K; \sll_\rho, \{h_{(1)}^\rho(\lambda), h_{(2)}^\rho\}, \mathfrak{o}} 
\right) \in \CC^*.
\end{equation}
It is an invariant of knots. 
Moreover, 
if $\rho_1$ and $\rho_2$ are two $\lambda$-regular representations 
which have the same character, 
then 
$\mathbb{T}^K_\lambda(\rho_1) = \mathbb{T}^K_\lambda(\rho_2)$. 
Thus the Reidemeister torsion $\mathbb{T}^K_\lambda$ defines a map on 
the set 
$X^{\mathrm{irr}}_\lambda(M_K) 
= \{\chi \in X^{\mathrm{irr}}(M_K) \; |\; 
\chi \text{ is } \lambda\text{-regular}\}$ of $\lambda$-regular characters.

\begin{remark}\label{rmkinv}
The Reidemeister torsion $\mathbb{T}^K_\lambda(\rho)$ defined 
in Equation~(\ref{Tordef}) is exactly the inverse of the one considered 
in Dubois~\cite{JDFibre}.
\end{remark}

\subsection{Review on 
the acyclic Reidemeister torsion for 
knot exteriors}
\label{section:acyclic_torsion}
We review the results of the Reidemeister torsion for 
acyclic local systems of knot exteriors in this section.
Let $K$ be a knot in a homology three sphere $M$
and $M_K$ its exterior.
\subsubsection*{The acyclic Reidemeister torsion 
of a knot exterior for
abelian representations}
Let $\psi_z$ be a homomorphism from $\pi_1(M_K)$ to $\C^*$
such that $\psi_z(\mu) = e^z$ 
where $z$ is a complex number and $\mu$ is the meridian of $K$.
We let $C_*(M_K; \C_{\psi_z})$ denote the following local system:
\[
  \C_{\psi_z} \otimes_{\Z[\pi_1(M_K)]} C_*(\widetilde M_K; \Z)
\]
where $\widetilde M_K$ is the universal cover of $M_K$
and $\C_{\psi_z}$ is a right $\Z[\pi_1(M_K)]$-module 
via the homomorphism $\psi$, i.e.,
$w \cdot \gamma = \psi_z(\gamma)^{-1} w$ 
for any $w \in \C$ and $\gamma \in \pi_1(M_K)$.

It is known that 
the torsion of $C_*(M_K; \C_{\psi_z})$ can be obtained from  
the normalized Alexander polynomial $\Delta_K(t)$ of $K$ as follows.

\begin{theorem}[Corollary 11.9 of Turaev~\cite{Turaev:2001}]
\label{thm:abelian_torsion}
If $z$ is a complex number such that $\Delta_K(e^z) \not = 0$,
then the complex $C_*(M_K; \C_{\psi_z})$ is acyclic and 
$\mathrm{Tor}(C_*(M_K; \C_{\psi_z}), \mathbf{c}^*_{\mathcal B})$ is equal to 
\[
\epsilon \cdot e^{nz/2} \frac{\Delta_K(e^z)}{e^{z/2}-e^{-z/2}}
\]
where $\epsilon \in \{\pm 1\}$, $n$ is some integer and 
$\mathcal B$ is a basis of the Lie algebra of $\C^*$, 
i.e., some non-zero element in $\C$.
\end{theorem}

We can regard the following function on $X^{ab}(M_K)$ as 
the Reidemeister torsion.
\[
  X^{ab}(M_K) \ni \chi_{\varphi_{z}} \mapsto
  \frac{\Delta_K(e^{2z})}{e^{z}-e^{-z}} \in \C.
\]
\subsubsection*{The acyclic Reidemeister torsion 
of a knot exterior for $\SL$-representations}
Let $\alpha$ be the abelianization homomorphism of $\pi_1(M_K)$
which send the meridian $\mu$ to $t$.
Let $\rho$ be an $\SL$-representation of $\pi_1(M_K)$.
We let $C_*(M_K; \C(t) \otimes \sll_\rho)$ denote the following local system:
\[
  (\C(t) \otimes \sll_{\rho}) 
    \otimes_{\Z[\pi_1(M_K)]} C_*(\widetilde M_K ; \Z)
\]
where $\widetilde M_K$ is the universal cover of $M_K$ and 
$\C(t) \otimes \sll_{\rho}$ is a right $\Z[\pi_1(M_K)]$-module
via the action $\alpha \otimes (Ad\circ \rho)$, i.e.,
$(f(t) \otimes v)\cdot \gamma 
= f(t)t^{\alpha(\gamma)} \otimes Ad_{\rho(\gamma)^{-1}}(v)$
for any $f(t)\in \C(t)$, $v \in \sll$ and $\gamma \in \pi_1(M_K)$.
For simplicity of notation, we let $\tsll_{\rho}$ stand for $\C(t) \otimes \sll_{\rho}$.

The following proposition holds for this chain complex.
\begin{proposition}[Proposition 3.1.1 in Yamaguchi~\cite{Yamaguchi}]
\label{prop:acyclicity}
If an $\SL$-representation $\rho$ is $\lambda$-regular, 
then 
$C_*(M_K; \tsll_{\rho})$
is acyclic.
\end{proposition}

\begin{theorem}[Kirk--Livingston~\cite{KL}, Kitano~\cite{Kitano}]\label{thm:Kitano_KL}
Let $\mathcal B$ be a basis of $\sll$.
If $C_*(M_K; \tsll_{\rho})$ is acyclic,
then the torsion ${\rm Tor}(C_*(M_K; \tsll_\rho),\mathbf{c}^*_{\mathcal B})$
coincides with the twisted Alexander invariant
of $\pi_1(M_K)$ and $Ad\circ \rho$.
\end{theorem}
The twisted Alexander invariant is given by using Fox differentials.
We will review it in the next section.

\section{The non-acyclic Reidemeister torsion at 
bifurcation points}
\label{section:torsion_at_bifurcation}
In this section, we will see that 
the limit of
the Reidemeister torsion ${\mathbb T}^{K}_{\lambda}$
is given by the differential coefficient of the acyclic Reidemeister torsion
$\Delta_K(e^{2z})/(e^z - e^{-z})$
at bifurcation points corresponding to simple roots of the Alexander
polynomial of $K$.
Here $\Delta_K(t)$ is normalized, i.e.,
$\Delta_K(t)=\Delta_K(t^{-1})$ and $\Delta_K(1)=1$.

\subsection{On the existence of a path of 
$\gamma$-regular characters}
\label{section:existence_limit}
We show that there exists a path of characters 
of $\gamma$-regular representations which converges to
a bifurcation point if the function $I_{\gamma}$ is not constant
on $X^{nab}(M_K)$ near the bifurcation point.

\begin{proposition}\label{prop:path_to_bifurcation}
Let $z_0$ be a complex number such that $e^{2z_0}$ is a simple root 
of the Alexander polynomial of $K$ and $\rho_{z_0}$ be a reducible 
non-abelian $\SL$-representation whose character is the same as that of 
the abelian representation $\varphi_{z_0}$.
Let $\gamma$ denote a simple closed curve in $\partial M_K$.
If the function $I_{\gamma}$ is not constant on the component of $X^{nab}(M_K)$ 
which contains the character $\chi_{\rho_{z_0}}$, 
then 
there exists a neighbourhood $V$ of $\chi_{\rho_{z_0}}$ such that
any point of $V$ except for at most finite points is $\gamma$-regular.
\end{proposition}

We prepare some notions to prove Proposition \ref{prop:path_to_bifurcation}.
Let $\rho$ be an irreducible $\SL$-representation of $\pi_1(M_K)$ such that 
$\rho(\pi_1(\partial M_K))$ contains a non-trivial hyperbolic element of $\SL$.
Let $\gamma$ be a simple closed curve in $\partial M_K$.
We can choose a neighbourhood $U$ of $\chi_{\rho}$ such that 
for any $\rho' \in \pi^{-1}(U)$, 
the image of the peripheral subgroup $\rho'(\pi_1(\partial M_K))$ also contains a non-trivial hyperbolic element.
We can define an analytic function $\alpha_{\gamma}$
on $U$ by the following equation:
\[
\rho'(\gamma) = 
A_{\rho'}
\left(
\begin{array}{cc}
e^{\alpha_{\gamma}} & 0 \\
 0 & e^{-\alpha_{\gamma}}
\end{array}
\right)
A^{-1}_{\rho'}
\]
where $A_{\rho'} \in \SL$ (for details, see Porti~\cite[Definition 3.19]{Porti:1997}).
Note that this function satisfies the following equation: 
\[
e^{2\alpha_{\gamma}(\chi)} - I_{\gamma}(\chi)e^{\alpha_{\gamma}(\chi)}+1=0.
\]

Proposition 3.26 in Porti~\cite{Porti:1997} gives a criterion about
the $\gamma$-regularity of $\rho$.

\begin{lemma}[Consequence of Proposition 3.26 in Porti~\cite{Porti:1997}]
\label{lemma:Porti_consequence}
Suppose that the dimension of the component containing $U$ is equal to $1$.
The irreducible representation $\rho$ is $\gamma$-regular if and only if
$\alpha_{\gamma}\circ \pi : \pi^{-1}(U) \subset R(\pi_1(M_K), \SL) \to \C$
is a submersion at $\rho$.
\end{lemma}

Proposition \ref{prop:path_to_bifurcation} follows
from Theorem \ref{thm:bifurcation_point} and 
Lemma \ref{lemma:Porti_consequence}.

\begin{proof}[Proof of Proposition \ref{prop:path_to_bifurcation}]
We let $X_0$ denote the component of $X^{nab}(M_K)$
which contains the bifurcation point $\chi_{\rho_{z_0}}$.
Theorem \ref{thm:bifurcation_point} implies that
the dimension of $X_0$ is equal to $1$.
Since $e^{2z_o}$ is a root of the Alexander polynomial of $K$,
$I_{\mu}(\chi_{\rho_{z_0}})$ is not equal to $\pm 2$.
In particular, $\rho_{z_0}(\mu)$ is a hyperbolic element in $\SL$.  
Thus the subgroup $\rho_{z_0}(\pi_1(\partial M_K))$ consists of  
hyperbolic elements.

By continuity, we can take a neighbourhood $U$ of $\chi_{\rho_{z_0}}$ in $X_0$
such that, for every $\chi \in U$,  
$I_{\mu}(\chi) \not = \pm 2$.
Let $V$ be a compact neighbourhood of $\chi_{\rho_{z_0}}$ in $U$.
Since $\alpha_{\gamma}$ is analytic and $I_{\gamma}$ is not constant in $V$, 
there exist only finite characters where the derivation of $\alpha_\gamma$ vanishes. 
Hence, by Lemma \ref{lemma:Porti_consequence}, 
there are only a finite number of characters in $V$ 
which are not $\gamma$-regular.
\end{proof}

\begin{cor}\label{cor:path_lambda}
If the function $I_{\lambda}$ is not constant near $\chi_{\rho_{z_0}}$
on $X^{nab}(M_K)$,
then there exists a path of $\lambda$-regular characters
which converges to $\chi_{\rho_{z_0}}$.
\end{cor}

\subsection{
Limits of the non-acyclic Reidemeister torsion
for knots at bifurcation points
}
\label{section:limit_torsion_formula}

If the function $I_{\lambda}$ is not constant near a bifurcation point
corresponding to a simple root of $\Delta_K(t)$, 
then there exists a path of $\lambda$-regular characters,
converging to the bifurcation point.
We can consider the limit of the Reidemeister torsion ${\mathbb T}^K_{\lambda}$
along this path.
This limit is obtained from the differential coefficient of
${\Delta_K(e^{2z})}/(e^{z}-e^{-z})$ as follows.
\begin{theorem}\label{thm:Main_theorem}
Let $z_0$ be a complex number
such that $e^{2z_0}$ is a simple root of 
the Alexander polynomial $\Delta_K(t)$ of K.
Let $\rho_{z_0}$ denote the reducible non-abelian $\SL$-representation
whose character is the same as one of $\varphi_{z_0}$.
If the function $I_{\lambda}$ is not constant 
near $\chi_{\rho_{z_0}}$ on $X^{nab}(M_K)$, 
then the limit of the Reidemeister torsion ${\mathbb T}^K_{\lambda}$
is expressed as 
\begin{equation}\label{eqn:main_equation}
\lim_{\chi_{\rho} \to \chi_{\rho_{z_0}}}
\mathbb T^{K}_{\lambda}(\rho)
= \varepsilon \cdot
  \left(
    \frac{1}{2}
    \left.
    \frac{d}{dz}
      \left(
        \frac{\Delta_{K}(e^{2z})}{e^z - e^{-z}}
      \right)
    \right|_{z=z_0}
  \right)^2
\end{equation}
where $\varepsilon \in \{\pm 1\}$.
\end{theorem}
The function ${\Delta_{K}(e^{2z})}/{(e^z - e^{-z})}$ is regarded as
the Reidemeister torsion for the abelian representation $\psi_z$
by Theorem \ref{thm:abelian_torsion}.
This relation shows that 
the Reidemeister torsion for the non-abelian representation $\rho_{z_0}$
is determined by the Reidemeister torsion for the abelian representation 
$\psi_{z_0}$.

\subsection{
Proof of Theorem \ref{thm:Main_theorem}
}
To prove this theorem, 
we describe the Reidemeister torsion $\mathbb T^K_{\lambda}(\chi_{\rho})$
as the differential coefficient of the sign--determined Reidemeister torsion
of $C_*(M_K; \tsll_{\rho})$ 
as follows. (see Theorem 3.1.2 in Yamaguchi~\cite{Yamaguchi}.)
\[
{\mathbb T}^K_{\lambda}(\chi_{\rho})
=
- \lim_{t \to 1} 
\frac{\T(M_K; \tsll_\rho, \mathfrak o)}{t-1}.
\]
where we write $\T(M_K; \tsll_\rho, \mathfrak o)$ instead of
$\mathrm{TOR}(M_K; \tsll_{\rho}, \emptyset, \mathfrak{o})$ for simplicity.
We want to know the following limit;
\[
\lim_{\rho \to \rho_{z_0}}
  \left(
    \lim_{t \to 1}
      \frac{\T(M_K; \tsll_\rho, \mathfrak o)}{t-1}
   \right).
\] 
Here we take the limit along a path of $\lambda$-regular representations,
converging to the reducible representation $\rho_{z_0}$.
We investigate the behavior of 
$\frac{\T(M_K; \tsll_\rho, \mathfrak o)}{t-1}$ 
at $\rho=\rho_{z_0}$ and $t=1$.
Since 
the numerator is regarded as the sign--determined twisted Alexander invariant
for $K$ and $Ad \circ \rho$ 
(for the details, see Kirk--Livingston~\cite{KL}, Kitano~\cite{Kitano} and Yamaguchi~\cite{Yamaguchi}), 
it is described more explicitly as follows.
Suppose that the group of $K$ has the following presentation:
\[
\pi_1(M_K)
=
\langle
x_1, \ldots, x_k \,|\, r_1, \ldots, r_{k-1}
\rangle.
\]
Since $\alpha:\pi_1(M_K) \to \Z$ is surjective,
by interchange columns if necessary,
we can assume that $\alpha(x_1) \not = 1$.
Then we have that $\det \Phi(x_1 - 1) \not = 1$.
By Proposition \ref{prop:acyclicity} and Theorem \ref{thm:Kitano_KL},
if a representation $\rho$ of $\pi_1(M_K)$ is $\lambda$-regular,
then the chain complex $C_*(M_K;\tsll_{\rho})$ is
acyclic and 
its torsion $\T(M_K;\tsll_{\rho}, \mathfrak o)$
is well--defined and given by  
\begin{equation}\label{torsion}
\tau_0 \cdot t^m
\frac{\det A^1_{K, Ad \circ \rho}}{\det \Phi(x_1-1)},
\end{equation}
where 
$m$ is some integer,
the symbol $\Phi$ stands for the tensor product homomorphism 
\[
  \alpha \otimes Ad \circ \rho: \Z[\pi_1(M_K)] \to M_3(\C[t, t^{-1}])
\]
with respect to a basis of $\sll$
and 
$A^1_{K, Ad\circ\rho}$ denotes
the following $3(k-1)\times 3(k-1)$ matrix over $\C[t, t^{-1}]$:
\[
A^1_{K, Ad\circ\rho}
=
\left(
\begin{array}{ccc}
\Phi(\frac{\partial r_1}{\partial x_2}) & \ldots & \Phi(\frac{\partial r_{k-1}}{\partial x_2}) \\
\vdots & \ddots & \vdots \\
\Phi(\frac{\partial r_1}{\partial x_k}) & \ldots & \Phi(\frac{\partial r_{k-1}}{\partial x_k})
\end{array}
\right).
\]
This rational function is the twisted Alexander invariant 
defined by Wada~\cite{wada}.
He has shown that 
the twisted Alexander invariant does not depend on the presentation of the group. (Theorem 1 in Wada~\cite{wada})
By the Euclidean algorithm, 
we can choose the following presentation for the knot group $\pi_1(M_K)$.
\begin{lemma}[Lemma 2.1 in Heusener--Porti--Su\'arez~\cite{HPS}]\label{lem:pres_group}
If necessary, we can replace the presentation of $\pi_1(M_K)$ by
$
\langle 
x'_1, \ldots, x'_k \,|\, r'_1, \ldots, r'_{k-1}
\rangle
$
such that $\alpha(x'_i) = t$ for all $i$.
\end{lemma}
\begin{remark}
The chosen presentation is not required to be a Wirtinger presentation 
in the case of a knot in $S^3$.
\end{remark}
Therefore we can assume  from the beginning that 
$\pi_1(M_K)$ has the presentation: 
$$\langle
x_1, \ldots, x_k \,|\, r_1, \ldots, r_{k-1}
\rangle$$
such that $\alpha(x_i) = t$ for all $i$.


The rational function $(\ref{torsion})$ is expressed as  
\[
\tau_0 \cdot
\frac{\det A^1_{K, Ad \circ \rho}}{\det \Phi(x_1-1)}
=
\frac{ \det A^1_{K, Ad \circ \rho} }{(t-1)(t^2 - \trace(\rho(x_1^2)) t +1)}.
\]
Therefore the torsion $\T(M_K, \tsll_{\rho}, \mathfrak o)/(t-1)$ is equal to
\[
\tau_0\cdot
\frac{\det A^1_{K, Ad \circ \rho}}{(t-1)^2 (t^2 - \trace(\rho(x_1^2)) t +1)}
\]
up to a factor $t^m$.
Since we suppose that $\rho$ is $\lambda$-regular, 
we know that $(t-1)^2$ divides $\det A^1_{K, Ad\circ \rho}$
(See Section 3.3 in Yamaguchi~\cite{Yamaguchi}).

Let $G_{\rho}(t)$ denote the rational function 
$(\det A^1_{K, Ad\circ \rho}) / (t-1)^2$.
We will consider the following two functions $t^2 - \trace(\rho(x_1^2)) t + 1$ 
and $G_{\rho}(t)$ at $\rho=\rho_{z_0}$ and $t=1$.


\begin{lemma}\label{lemma:denominator}
The function $t^2 - \trace(\rho_{z_0}(x_1^2)) t + 1$ is smooth and non-zero 
at $\rho=\rho_{z_0}$ and $t=1$.
\end{lemma}
\begin{proof}[Proof of Lemma \ref{lemma:denominator}]
The function $t^2 - \trace(\rho(x_1^2)) t + 1$ depends on $\rho$ 
smoothly.
We look for the value of $t^2 - \trace(\rho(x_1^2)) t + 1$ at 
$\rho = \rho_{z_0}$.
By the assumption that $\rho_{z_0}$ has the same character as $\varphi_{z_0}$, 
we have that 
$
\trace(\rho_{z_0}(x_1^2)) = e^{2z_0} + e^{-2z_0}.
$ 
Since $e^{2z_0}$ is a simple root of the Alexander polynomial $\Delta_K(t)$ of $K$
and $\Delta_K(1)=1$,
the complex number $e^{2z_0}$ is not equal to $1$.
Hence if we substitute $t=1$ into 
the polynomial $t^2 - \trace(\rho_{z_0}(x_1^2)) t + 1$, 
then its value $2-(e^{2z_0}+e^{-2z_0})$ is not zero.
\end{proof}

The following proposition plays an important role 
when we consider the function $G_{\rho}(t)$ and 
prove Theorem \ref{thm:Main_theorem}.

\begin{proposition}\label{prop:t_Alexander_bifurcation}
The chain complex $C_*(M_K; \tsll_{\rho_{z_0}})$
is acyclic.
Moreover the Reidemeister torsion 
$\T(M_K; \tsll_{\rho_{z_0}}, \mathfrak o)$ 
is given by 
\begin{equation}\label{eqn:t_Alexander_bifurcation}
\tau_0
\cdot
\epsilon t^m
\cdot
\frac{\Delta_K(t)\Delta_K(te^{2z_0})\Delta_K(te^{-2z_0})}
{(t-1)(t^2 - \trace(\rho_{z_0}(x_1^2)) t + 1)}
\end{equation}
where 
$\epsilon \in \{\pm 1\}$, $m \in \Z$ and 
$\Delta_K(t)$ is the normalized Alexander polynomial of $K$.
\end{proposition}

\begin{proof}[Proof of Proposition \ref{prop:t_Alexander_bifurcation}]
It is enough to prove 
the following claims:
\begin{itemize}
\item $\det(\Phi(x_1)-1)$ is not zero;
\item $\det A^{1}_{K, Ad\circ\rho_{z_0}}$ is expressed by using the product of the three Alexander polynomials which appear in the numerator of the fraction in Eq.\,$(\ref{eqn:t_Alexander_bifurcation})$;
\item $C_*(M_K;\tsll_{\rho_{z_0}})$ is acyclic and its Reidemeister torsion is given as above.
\end{itemize}

We have seen that $\det( \Phi(x_1 -1) )$ is not zero.
We consider $\det A^{1}_{K, Ad\circ\rho_{z_0}}$.
Since the $\SL$-representation $\rho_{z_0}$ has the same character as 
$\varphi_{z_0}$ and $\alpha(x_i) = \alpha(\mu)$ for all $i$, 
we have that 
\[
\trace(\rho_{z_0}(x_1))
  =\cdots
  =\trace(\rho_{z_0}(x_k))
  =\trace(\rho_{z_0}(\mu)).
\]
Furthermore $\rho_{z_0}$ is reducible, then we can assume that 
\[
\rho_{z_0}(x_i)
=
\left(
\begin{array}{cc}
e^{z_0} & \alpha_i \\
0 & e^{-z_0} 
\end{array}
\right)
\]
by taking conjugation, where $\alpha_i$ is a complex number 
(Remark \ref{rem:maximal_sub}).
We take an ordered basis $\{E, H, F\}$ of $\sll$ as follows:
\[
E=
\left(
\begin{array}{cc}
0 & 1 \\
0 & 0
\end{array}
\right), 
H=
\left(
\begin{array}{cc}
1 & 0 \\
0 & -1
\end{array}
\right), 
F=
\left(
\begin{array}{cc}
0 & 0 \\
1 & 0
\end{array}
\right).
\]
Under this basis, 
for each $x_i$, 
the representation matrix of $Ad(\rho_{z_0}(x_i) )$
is given by
\[
Ad(\rho_{z_0}(x_i)^{-1})
=
\left(
\begin{array}{ccc}
e^{-2z_0} & 2\alpha_i e^{-z_0} & -\alpha_i^2 \\
0 & 1 & -\alpha_i e^{z_0} \\
0 & 0 & e^{2z_0} 
\end{array}
\right).
\]

Note that
each $\Phi(\frac{\partial r_i}{\partial x_j })$
is an upper triangular matrix for any $i$ and $j$.

We express the matrix   
$\Phi( \frac{\partial r_i}{\partial x_j} )\, 
(1 \leq i \leq k-1, 2 \leq j \leq k)$
by using the following matrix:
\[
\left(
\begin{array}{ccc}
a_{ij} & * & * \\
 0 & b_{ij} & * \\
 0 & 0 & c_{ij}
\end{array}
\right).
\]

\begin{claim}\label{claim:1st_torsion}
\[
\det A^1_{K, Ad \circ \rho_{z_0}} 
=
\epsilon t^m
\Delta_K(t) \Delta_K(te^{2z_0}) \Delta_K(te^{-2z_0})
\]
where $\Delta_K(t)$ is the normalized Alexander polynomial of $K$, 
$\epsilon \in \{\pm 1\}$ and $m \in \Z$.
\end{claim}

\begin{proofclaim}[Proof of Claim \ref{claim:1st_torsion}]
\begin{align*}
\det A^1_{K, Ad \circ \rho_{z_0}}
&=
\left|
\begin{array}{ccccccc}
a_{12} & * & *  & a_{22} & * & * & \cdots  \\
  & b_{12} & *  &   & b_{22} & * & \vdots \\
  &   & c_{12}  &   &   & c_{22} & \cdots \\
a_{13} & * & * & a_{23} & * & * & \cdots  \\
  & b_{13} & * &   & b_{23} & * & \vdots \\
  &   & c_{13} &   &   & c_{23} & \cdots  \\
  & \vdots & & & \vdots     &    & \ddots  
\end{array} 
\right|  \\
&=
\left|
\begin{array}{ccc}
A & * & * \\
  & B & * \\
  &   & C 
\end{array}
\right|.
\end{align*}
Here $A, B$ and $C$ 
respectively denote the small matrices 
$(a_{ij})_{i, j}, (b_{ij})_{i,  j}$
and $(c_{ij})_{i, j}$\, $(1\leq i \leq k-1, 2\leq j \leq k)$.

From the equation $\alpha(x_1) = t$ and 
the calculation of the Alexander polynomial using Fox differentials 
(Chapter $9$ in Burde and Zieschang~\cite{Burdeknots}), 
we can see that 
there exist some integer $n'$ and $\epsilon \in \{\pm 1\}$ such that
\begin{align*}
\det A &= \epsilon (e^{-2z_0}t)^{n'}  \Delta_K(te^{-2z_0}), \\
\det B &= \epsilon t^{n'}  \Delta_K(t), \\
\det C &= \epsilon (e^{2z_0}t)^{n'}  \Delta_K(te^{2z_0}).
\end{align*}
Therefore we have that
\[
\det A^1_{K, Ad \circ \rho_{z_0}}
=
\epsilon t^{3n'} 
\Delta_K(t)\Delta_K(te^{2z_0})\Delta_K(te^{-2z_0}).
\]
\hfill (Claim \ref{claim:1st_torsion})
\end{proofclaim}

\noindent
Hence $C_*(M_K; \tsll_{\rho_{z_0}})$ is acyclic. 
Furthermore we can see that there exists some integer $m$ 
such that the sign--determined Reidemeister torsion of 
$C_*(M_K; \tsll_{\rho_{z_0}})$ is expressed as  
\[
\mathcal{T}(M_K;\tsll_{\rho_{z_0}}, \mathfrak o)
=
\tau_0
\cdot
\epsilon t^m
\cdot
\frac{\Delta_K(t)\Delta_K(te^{2z_0})\Delta_K(te^{-2z_0}) }
     {(t-1)(t^2 - \trace(\rho_{z_0}(x_1^2)) t + 1)}.
\]
\end{proof}

Now we consider the function $G_{\rho}(t)$ at $\rho=\rho_{z_0}$ and $t=1$.
\begin{lemma}\label{lemma:numerator}
The rational function 
$G_{\rho}(t)$ is smooth and non-zero 
at $\rho=\rho_{z_0}$ and $t=1$.
\end{lemma}
\begin{proof}[Proof of Lemma \ref{lemma:numerator}]
Since $e^{2z_0}$ is a simple root of $\Delta_K(t)$ and
$\Delta_K(t)$ is symmetric for $t$, i.e., $\Delta_K(t)=\Delta_K(t^{-1})$, 
the complex number $e^{-2z_0}$ is also a simple root of $\Delta_K(t)$.
Hence the numerator $\Delta_K(t)\Delta_K(te^{2z_0})\Delta_K(te^{-2z_0})$ 
of Eq.\,$(\ref{eqn:t_Alexander_bifurcation})$ has the second order zero at $t=1$.
Therefore 
the function 
$\det A^1_{K, Ad\circ \rho}$
can be divided by $(t-1)^2$ at $\rho=\rho_{z_0}$.
We can define $G_{\rho_{z_0}}(t) \in \C[t, t^{-1}]$.
Hence the function $G_{\rho}(t)$ changes smoothly to $G_{\rho_{z_0}}(t)$
and 
there exists non-zero limit of $G_{\rho}(t)$ at $\rho=\rho_{z_0}$ and $t=1$.
\end{proof}

Now, we are ready to calculate the limit of the Reidemeister torsion ${\mathbb T}^K_{\lambda}$
by using Proposition \ref{prop:t_Alexander_bifurcation}.

\begin{proof}[Proof of Themorem \ref{thm:Main_theorem}]
By Lemmas \ref{lemma:denominator} and \ref{lemma:numerator}, 
we see that 
the limit of the rational function 
$\T(M_K; \tsll_{\rho}, \mathfrak o)/(t-1)$ 
at $\rho=\rho_{z_0}$ and $t=1$ exists.
Moreover when we express the rational function 
$\T(M_K; \tsll_{\rho}, \mathfrak o)/(t-1)$ 
as $G_{\rho}(t)/(t^2 - \trace(\rho(x_1^2))+1)$,
both of the numerator $G_\rho(t)$ and the denominator $t^2 - \trace(\rho(x_1^2))+1$ are smooth and non-zero near $\rho=\rho_{z_0}$ and $t=1$.
Hence we can change the order of taking limits.
By interchanging the limit of $t$ and that of $\rho$ and 
by Proposition \ref{prop:t_Alexander_bifurcation},
the limit of ${\mathbb T}^K_{\lambda}$ is calculated as follows.
\begin{align*}
\lim_{\rho \to \rho_{z_0}} {\mathbb T}^K_{\lambda}(\rho)
&=
-\lim_{\rho \to \rho_{z_0}}
\left(
\lim_{t\to 1}
\frac{\mathcal{T}(M_K;\tsll_{\rho}, \mathfrak o)}
     {t-1}
\right)\\
&=
-\lim_{t \to 1}
\frac{\mathcal{T}(M_K;\tsll_{\rho_{z_0}}, \mathfrak o)}
{t-1}\\
&=
 \lim_{t \to 1}\left\{
\frac{\Delta_K(te^{2z_0})\Delta_K(te^{-2z_0})}{(t-1)^2}
\cdot
\frac{-\epsilon \tau_0 t^m \Delta_K(t)}{t^2 - \trace(\rho_{z_0}(x_1^2)) t +1}
\right\}
\end{align*}
where $\epsilon \in \{\pm 1\}$ and $m \in \Z$.

Since $\Delta_K(1)=1$ and $e^{2z_0}$ and $e^{-2z_0}$ are simple roots of 
$\Delta_K(t)$, we have
\begin{align}
& \lim_{t \to 1}\left\{
\frac{\Delta_K(te^{2z_0}) \Delta_K(te^{-2z_0}) }{(t-1)^2}
\cdot
\frac{-\epsilon \tau_0 t^m \Delta_K(t)}{t^2 - \trace(\rho_{z_0}(x_1^2)) t +1} 
\right\}
\nonumber \\
&=
\frac{-\epsilon \tau_0}{2 - (e^{2z_0} + e^{-2z_0})}
\cdot
\lim_{t \to 1}
\left\{
\frac{\Delta_K( te^{2z_0} ) }{t-1} \cdot \frac{\Delta_K( te^{-2z_0} ) }{t-1} 
\right\}
\nonumber \\
&=
\frac{-\epsilon \tau_0 \Delta'_K(e^{2z_0}) \Delta'_K(e^{-2z_0})}
     {2 - (e^{2z_0} + e^{-2z_0})}. \label{eqn:limit_bifurcation}
\end{align}

It follows from the symmetry of $\Delta_K(t)$ that 
\begin{equation}\label{eqn:diff_Alex_poly}
\Delta'_K(e^{-2z_0}) = - \Delta'_K(e^{2z_0}) e^{4z_0}.
\end{equation}
If we substitute Equation $(\ref{eqn:diff_Alex_poly})$ 
into Equation $(\ref{eqn:limit_bifurcation})$,
then we obtain: 
\begin{equation*}
\lim_{\rho \to \rho_{z_0}}
\mathbb T^K_{\lambda} (\rho)
=
-\tau_0
\cdot
\epsilon
\cdot
\frac{ ( \Delta'_K(e^{2z_0}) e^{2z_0} )^2 }{(e^{2z_0}+e^{-2z_0})-2}.
\end{equation*}


On the other hand, 
the right hand side of Equation $(\ref{eqn:main_equation})$ 
is given by a direct calculation:
\begin{align*}
\left(
\frac{1}{2}
\left.
\frac{d}{dz}
\left(
\frac{\Delta_K(e^{2z})}{e^z-e^{-z}}
\right)
\right|_{z=z_0}
\right)^2
&=
\left(
\frac{\Delta'_K(e^{2z_0}) e^{2z_0}}{e^{z_0}-e^{-z_0}}
\right)^2\\
&=
\frac{(\Delta'_K(e^{2z_0}) e^{2z_0})^2}{e^{2z_0}+e^{-2z_0} -2 }.
\end{align*}
Therefore we have
\[
\lim_{\rho \to \rho_{z_0}}
\mathbb T^K_{\lambda} (\rho)
=  \varepsilon \cdot
\left(
  \frac{1}{2}
  \left.
    \frac{d}{dz}
    \left(
      \frac{\Delta_K(e^{2z})}{e^z-e^{-z}}
    \right)
  \right|_{z=z_0}
\right)^2
\]
where $\varepsilon = - \tau_0 \cdot \epsilon$,
which completes the proof.
\end{proof}

\section{
On the reference generators of 
the $\sll_\rho$-homology groups
at a bifurcation point}
\label{section:reference_generators}

We consider 
the reference generators of $H_*(M_K; \sll_\rho)$ in this section.
By Proposition \ref{prop:path_to_bifurcation},
the reference generators
$\{h^{\rho}_{(1)}(\lambda), h^{\rho}_{(2)}\}$ of $H_*(M_K; \sll_\rho)$
exist for any irreducible representation $\rho$
sufficiently near the reducible non-abelian representation $\rho_{z_0}$
when $I_\lambda$ is not constant on the component of $X^{nab}(M_K)$, 
containing the bifurcation point $\chi_{\rho_{z_0}}$.
Here the representation $\rho_{z_0}$ corresponds to 
a simple root of the Alexander polynomial of
$K$.
We will show that they can be extended to
the generator of $H_*(M_K; \sll_{\rho_{z_0}})$.
In this section, we assume that 
the regular function $I_\lambda$ is not constant
on the component containing the bifurcation point $\chi_{\rho_{z_0}}$
in $X^{nab}(M_K)$. 
\subsection{
On the generator of the second
$\sll_\rho$-twisted homology group 
at a bifurcation point.
}

From the following results 
of Heusener, Porti and Su\'arez~\cite{HPS} 
we know 
the dimensions of $H_*(M_K; \sll_{\rho_{z_0}})$ 
and the basis of $H_2(M_K; \sll_{\rho_{z_0}})$.

\begin{lemma}[Lemma 4.1 of Heusener--Porti--Su\'arez~\cite{HPS}]
Let $M$ be a connected, compact, orientable, irreducible 
$3$-manifold such that 
$\partial M$ is a torus and the first Betti number is one. 

Let $\rho : \pi_1(M) \to \SL$ be a non-abelian representation 
such that 
$\rho(\pi_1(\partial M))$ contains a non-parabolic element.
We let $i$ denote the inclusion $\partial M \hookrightarrow M$
and 
$Z^1(M; \sll_\rho)$ denote
the set of twisted cocycles of $M$ with coefficients in $\sll_\rho$.
If $\dim_{\C} Z^1(M; \sll_\rho)=4$, 
then we have an injection 
\[
i^* : H^1(M; \sll_\rho) \to H^1(\partial M; \sll_\rho)
\] 
and 
an isomorphism 
$i^* : H^2(M; \sll_\rho) \to H^2(\partial M; \sll_\rho)$.
\end{lemma}

We apply this lemma to the knot exterior $M_K$ and $\rho_{z_0}$.
It follows from Proposition 4.4 of Heusener--Porti--Su\'arez~\cite{HPS} 
that $\dim_{\C}Z^1(M_K; \sll_{\rho_{z_0}})=4$.
Since $\trace(\rho_{z_0}(\mu))=e^{z_0}+e^{-z_0} \not = \pm 2$, 
we see that 
$\rho_{z_0}(\pi_1(\partial M_K))$ contains a non-parabolic element.

Therefore we have that:
\begin{itemize} 
\item
  $H_0(M_K; \sll_{\rho_{z_0}}) = 0$; 
\item
  the induced homomorphism
  $i_* : H_1(\partial M_K; \sll_{\rho_{z_0}}) \to H_1(M_K; \sll_{\rho_{z_0}})$
  is surjective;
\item
  the induced homomorphism
  $i_* : H_2(\partial M_K; \sll_{\rho_{z_0}}) \to H_2(M_K; \sll_{\rho_{z_0}})$
  is an isomorphism.
\end{itemize}

Note that 
$\dim_{\C}H_2(\partial M_K; \sll_{\rho_{z_0}})$ is equal to $1$
since the restriction of $\rho_{z_0}$ to $\pi_1(\partial M_K)$ is non-trivial. 

\begin{proposition}\label{prop:base_second_homology}
The chain $i_*(P^{\rho_{z_0}} \otimes \widetilde{\partial  M_K})$
determines a basis of the homology group $H_2(M_K; \sll_{\rho_{z_0}})$.
Here $P^{\rho_{z_0}}$ is a vector in $\sll$ 
such that 
$Ad_{\rho_{z_0}(\gamma)}(P^{\rho_{z_0}}) = P^{\rho_{z_0}}$ for all 
$\gamma \in \pi_1(\partial M_K)$. 
\end{proposition}
\begin{proof}[Proof of Proposition \ref{prop:base_second_homology}]
By calculations, 
we see that 
$P^{\rho_{z_0}} \otimes \widetilde{\partial  M_K}$ is a cycle in 
$C_2(\partial M_K; \sll_{\rho_{z_0}})$ and
it determines a non-zero element of $H_2(M_K; \sll_{\rho_{z_0}})$
(see Porti~\cite[Proposition 3.18]{Porti:1997}).

Since $[P^{\rho_{z_0}} \otimes \widetilde{\partial M_K}]$ is a generator of 
$H_2(\partial M_K; \sll_{\rho_{z_0}})$, 
we can take $i_* ([P^{\rho_{z_0}} \otimes \widetilde{\partial M_K}])$ as 
a generator of $H_2(M_K; \sll_{\rho_{z_0}})$.
\end{proof}

\subsection{ 
On the generator of 
the first $\sll_\rho$-twisted homology group 
at a bifurcation point.
}

As $\partial M_K$ is a two--dimensional torus,
it follows from calculations that $H_1(\partial M_K; \sll_{\rho_{z_0} })$ 
is generated by 
$[P^{\rho_{z_0} } \otimes \widetilde \mu]$ and 
$[P^{\rho_{z_0}} \otimes \widetilde \lambda]$
(see Porti~\cite[Proposition 3.18]{Porti:1997}).
The problem lies in whether 
$i_* ([P^{\rho_{z_0} } \otimes \widetilde \lambda])$ is zero or not in 
$H_1(M_K; \sll_{\rho_{z_0} })$.
We shall show that 
$i_* ([P^{\rho_{z_0} } \otimes \widetilde \lambda])$ is a non-zero class in 
$H_1(M_K; \sll_{\rho_{z_0} } )$.
This follows from the fact 
that the limit of the Reidemeister torsion ${\mathbb T}^K_{\lambda}$ is not zero.
Together with Proposition \ref{prop:base_second_homology},
the following proposition holds.

\begin{proposition}\label{prop:extend_ref_gen}
Let $z_0$ be a complex number 
such that $e^{2z_0}$ is a simple root of the Alexander polynomial of $K$.
Let $\rho_{z_0}$ denote
the reducible non-abelian $\SL$-representation 
whose character is the same as
that of $\varphi_{z_0}$.
If $I_{\lambda}$ is not constant near the bifurcation point $\chi_{\rho_{z_o}}$,
then the reference generators $h_{(1)}^{\rho}(\lambda)$ and $h_{(2)}^{\rho}$
can be extended in $H_*(M_K; \sll_{\rho_{z_0}})$.
\end{proposition}

\begin{proof}[Proof of Proposition \ref{prop:extend_ref_gen}]
It is enough to show that 
$i_*([P^{\rho_{z_0}} \otimes \widetilde \lambda])$ is a non-zero class 
in $H_1(M_K; \sll_{\rho_{z_0}})$.
To this purpose, suppose that 
$i_*([ P^{\rho_{z_0}} \otimes \widetilde \lambda])$ is zero in 
$H_1(M_K; \sll_{\rho_{z_0}})$.

By Theorem \ref{thm:bifurcation_point},
it follows that $\dim_{\C}H_1(M_K; \sll_{\rho_{z_0}})=1$ 
and $\rho_{z_0}$ is a smooth point in 
the $\SL$-representation variety of the knot group.
By Corollary \ref{cor:path_lambda}, there exists a path 
$\{\rho_s \,|\, s \in \C, |s| < \epsilon \}$ of $\SL$-representations
such that $\rho_0= \rho_{z_0}$ and $\rho_s$ is $\lambda$-regular at $s \not =0$.
Here $\epsilon$ is a small positive real number. 
The cohomology group $H^1(M_K;\sll_{\rho_s})$ is isomorphic to 
the Zariski tangent space of $X(M_K)$ at $\chi_{\rho_s}$.
We can take a smooth family of generators 
$\{ \xi_s \}_s$ of $H^1(M_K; \sll_{\rho_s})$
associated with the path $\{ \rho_s \}_s$.
Using the Kronecker pairing between 
the homology group $H_1(M_K; \sll_{\rho_s})$ and 
the cohomology group $H^1(M_K;\sll_{\rho_s})$,
we have a family $\{\sigma_s\}_s$ of generators of 
$H_1(M_K; \sll_{\rho_s})$
such that the Kronecker pairing for $\sigma_s$ and $\xi_s$
does not vanish for each $s \in \C, |s| < \epsilon$.

We define 
a non-zero complex number ${\mathbb T}^K_{\sigma}(\rho_s)$ for each $s$ to be
\[
{\mathbb T}^K_{\sigma}(\rho_s)
=
\mathrm{TOR}(M_K; \sll_{\rho_s}, \{\sigma_s, h^{\rho_s}_{(2)}\}, \mathfrak o).
\]
where 
$h^{\rho_s}_{(2)}$ is the reference generator 
of $H_2(M_K; \sll_{\rho_s})$.
This function depends on $s$ smoothly.

\begin{claim}\label{claim:calculation_base_change}
Let $c_s$ denote the ratio between $h^{\rho_s}_{(1)}(\lambda)$ and 
$\sigma_s$, i.e., $h^{\rho_s}_{(1)}(\lambda) = c_s \cdot \sigma_s$.
Then the following equation holds at $s \not = 0$:
\[
{\mathbb T}^K_{\lambda}(\rho_s)
=
c_s
\cdot
{\mathbb T}^K_{\sigma}(\rho_s).
\]
\end{claim}

\begin{proofclaim}[Proof of Claim \ref{claim:calculation_base_change}]
This follows from the base change formula for the Reidemeister torsion
(see Dubois~\cite[Formula (5)]{JDFourier} and Porti~\cite[Proposition 0.2]{Porti:1997}).
\begin{align*}
{\mathbb T}^K_{\lambda}(\rho_s) 
&=
\mathrm{TOR}(M_K; \sll_{\rho_s}, 
\{h^{\rho_s}_{(1)}(\lambda), h^{\rho_s}_{(2)}\}, \mathfrak o)\\
&=
\mathrm{TOR}(M_K; \sll_{\rho_s}, \{\sigma_s, h^{\rho_s}_{(2)}\}, \mathfrak o)
\cdot
[h^{\rho_s}_{(1)}(\lambda) / \sigma_s] \\
&=
{\mathbb T}^K_{\sigma}(\rho_s) \cdot c_s.
\end{align*}
\hfill (Claim \ref{claim:calculation_base_change})
\end{proofclaim}

The function ${\mathbb T}^K_{\lambda}(\rho_s)$ also depends on $s$ smoothly. 
We have known from Theorem \ref{thm:Main_theorem} that 
there exists the non-zero limit of 
${\mathbb T}^K_{\lambda}(\rho_s)$ taking limit $s$ to $0$.

On the other hand, 
the limit of $c_s$ at $s=0$ is zero by the assumption.
The function ${\mathbb T}^K_{\sigma}(\rho_s)$  
does not have a pole at $s=0$ 
by the construction.
Hence
if we take a limit of $s$ to $0$, 
the function $c_s \cdot {\mathbb T}^K_{\sigma}(\rho_s)$ must be zero.
This is a contradiction.
Therefore 
$i_*([P^{\rho_{z_0}} \otimes \widetilde \lambda])$ determines 
a non-zero class in $H_1(M_K; \sll_{\rho_{z_0}})$.
\end{proof}




\end{document}